# A Novel and Efficient Algorithm for Scanning All Minimal Cutsets of a Graph

**Ahmad R. Sharafat**

and

**Omid R. Ma'rouzi**

Department of Electrical Engineering
Tarbiat Modarres University, Tehran, Iran.

____________________________________________________________________

We propose a novel algorithm for enumerating and listing all minimal cutsets of a given graph. It is known that this problem is NP-hard. We use connectivity properties of a given graph to develop an algorithm with reduced complexity for finding all its cutsets. We use breadth first search (BFS) method in conjunction with edge contraction to develop the algorithm. We introduce the concepts of a pivot vertex and absorbable clusters and use them to develop an enhanced recursive contraction algorithm. The complexity of the proposed algorithm is proportionate to the number of cutsets. We present simulation results to compare the performance of our proposed algorithm with those of existing methods.

____________________________________________________________________

## I. INTRODUCTION

SCANNING all minimal cutsets is an important issue in many applications, such as evaluating the reliability of networks [Rai 1982; Colbourn 1987; Fard and Lee 1999] or calculating the maximum flow through a network [Goldberg 1988; Ahuja and Magnati 1993]. We have also used it to find a lower bound for the routing and wavelength assignment (RWA) problem in all-optical wide area networks [Sharafat and Ma'rouzi 2002]. A closely related problem is to find all $(s,t)$-mincuts of a graph in which $(s,t)$ is a vertex pair. This problem is equivalent to finding all minimal cutsets that separate vertices in a vertex pair $(s,t)$ from each other in the graph [Provan and Shier 1996]. It is

This research was supported in part by Tarbiat Modarres University.
Authors' addresses: Ahmad R. Sharafat, Department of Electrical Engineering, Tarbiat Modarres University, P.O. Box: 14155-4838, Tehran, Iran, e-mail: *sharafat@isc.iranet.net*; Omid R. Ma'rouzi, Department of Electrical Engineering, Tarbiat Modarres University, P.O. Box: 14155-4838, Tehran, Iran, e-mail: *ormzi@yahoo.com*

also extensively used for calculating the reliability of a 2-terminal network [Shier 1988; Ball 1995].

Provan and Ball [1983] proved that scanning all cutsets of a given graph is an NP-hard problem. The basic conventional approach for scanning all cutsets of a given graph uses state-space enumeration. The method proposed by Ahmad [1988] to enumerate the minimal cutsets of acyclic directed graphs, and its extension to all terminal undirected graphs proposed by Hongfen [1995] are two such examples. Although it is conceptually simple, the state-space approach is impractical because the size of the state space grows exponentially with the number of nodes. Certain improvements to this approach focus directly on topological properties of the graph to substantially reduce the size of state-space. Using this concept, Tsukiyama *et al.* [1980] proposed an algorithm for finding the (*s,t*)-mincuts of a given graph in linear time per cutset. This recursive algorithm enumerates all connected *s-t* separating cutsets by considering only certain extensions (1-point extensions) of the separating cutsets.

Using Tsukiyama algorithm, Whited *et al.* [1990] proposed an improved algorithm for enumerating (*s,t*)-mincuts of planar graphs. An efficient method for enumerating all (*s,t*)-mincuts in a directed graph is the paradigm of Provan and Shier [1996]. It is based on a pivotal decomposition on a vertex.

A fast random approach to enumerate all cutsets within a multiplicative factor of the minimum cutset was proposed by Karager [2000]. This algorithm is based on iterative contraction of edges in a graph and lists all cutsets of a graph while the relative error can be reduced arbitrarily by increasing the number of trials in the algorithm.

In this paper, we modify the method proposed by Tsukiyama *et al.* [1980] by using the concept of iterative contraction [Karger 2000] and BFS ordering of vertices to develop a novel recursive contraction algorithm for scanning (enumerating and listing) all minimal cutsets of a given graph. Also, we introduce the concepts of a pivot vertex and absorbable and inabsorbable clusters, and use them to develop an enhanced recursive contraction algorithm. We show that the complexity of this enhanced recursive contraction algorithm is pseudo-polynomial and is lower than the existing methods proposed by Ahmad [1988] and by Tsukiyama *et al.* [1980].

In Section II we present our notations and introduce some basic concepts. In Section III we develop a recursive contraction algorithm for finding all possible cutsets of a given graph in linear time per cutset. In Section IV we present the results of applying the proposed algorithm to different sample graphs, and finally in Section V we present the conclusion and a summary of results.

## II. PRELIMINARIES

An undirected graph $G=(V,E)$ consists of a set $V$ of vertices and a set $E$ of edges whose elements are unordered pairs of vertices. The edge $e=(u,v) \in E$ is said to be incident with vertices $u$ and $v$, where $u$ and $v$ are the end points of $e$. These two vertices are called adjacent. The set of vertices adjacent to $v$ is written as $A(v)$, and the degree of $v$ is the number of vertices adjacent to $v$ and is denoted as $|A(v)|$. Throughout, we will reserve $n$ for $|V|$ and $m$ for $|E|$. In a graph $G=(V,E)$, a partition $(X,X')$ is defined as the two proper disjoint subsets of $V$. The complement of $X \subseteq V$ is denoted as $X' = V-X$. The open neighborhood of $X$ is defined as $\Gamma(X)=\{v \in X' \mid (u,v) \in E \text{ for some } u \in X \}$. The induced subgraph $\langle X \rangle$ is the graph $H = (X,F)$ where $F=\{(u,v) \in E \mid u,v \in X\}$.

An alternating sequence of distinct adjacent vertices and their incident edges, i.e., $v_0$, $(v_0, v_1)$, $v_1$, ... , $v_{k-1}$, $(v_{k-1}, v_k)$, $v_k$ is called a *u-v* path when $v_0= u$ and $v_k = v$. If a *u-v* path exists in a graph $G$ between all vertices $u$ and $v$, then $G$ is connected. Otherwise, $G$ decomposes into a number of connected subgraphs referred to as components of G. A graph *G/{u}* obtained by deleting a vertex *u* is defined as a subgraph of *G* induced by *V-{u}*, i.e., *G/{u}* = ⟨*V-{u}*⟩.

**Definition 1. Cutset:** For a given graph *G= (V,E),* a subset of edges $C \tilde{I} E$ is a minimal cutset if and only if deleting all edges in *C* would divide *G* into two connected components. An isolated vertex is considered as a component. A minimal cutset divides the vertices of *G* into two disjoint subsets *X* and *X'*, each of which induce a connected subgraph. We denote a cutset as ⟨*X,X'*⟩. For convenience we use cutset instead of minimal cutset.

**Definition 2. Exterior and Interior Vertices:** For a connected graph *G*, a vertex *v* is called an interior vertex if the graph *G/{v}* = ⟨*V - {v}*⟩ is not connected, otherwise it is called an exterior vertex. In Fig. 1, vertex 2 is an interior vertex and vertex 3 is an exterior vertex.

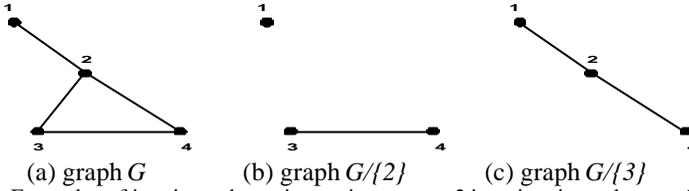

(a) graph *G*      (b) graph *G/{2}*      (c) graph *G/{3}*

**Fig.1** - Examples of interior and exterior vertices; vertex 2 is an interior and vertex 3 is an exterior vertex

**Definition 3. Monolithic and multi-segment graphs:** A monolithic graph does not have any interior vertex, and a multi-segment graph has at least one interior vertex. If *v* is an interior vertex of a multi-segment graph *G*, its deletion results in multiple components $C_i = \langle H_i \rangle$. We call $S_i = \langle H_i \tilde{E}\{v\} \rangle$ a segment. Each segment $S_i$ can also be a multi-segment graph. A monolithic segment is also called an elementary segment. Fig. 2 shows a multi-segment graph and its elementary segments.

**Observation 1.** A multi-segment graph can be decomposed into a number of elementary segments. Each cutset of an elementary segment is a cutset of a multi-segment graph. So the set of all cutsets of a multi-segment graph is the union of the sets of all cutsets of its elementary segments. Without any loss of generality we assume that graphs are monolithic.

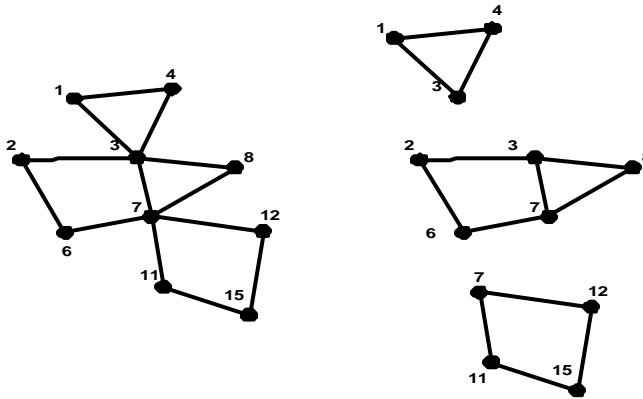

**Fig. 2** - (a) multi-segment graph *G*,      (b) elementary segments of *G*

**Definition 4. Edge Contraction:** A graph denoted as *G/uv* is made by the contraction of an edge *uv* in the Graph *G* in the following manner. Delete vertices *u* and *v* in *G* and replace them by a new contracted vertex *g*. Then remove all edges incident to both *u* and *v* (i.e., *uv*). For each edge incident to one of the vertices *u* or *v* (i.e., *uw* or *vw*), there is an incident edge between *g* and the other vertex of the incident edge (i.e., *gw*) in *G/uv*. We extend this definition for an edge set *F* by sequentially applying the contraction to all edges of *F* (the order of contractions is irrelevant). The resulting graph is denoted as *G/F*. An example of edge contraction is shown in Fig. 3.

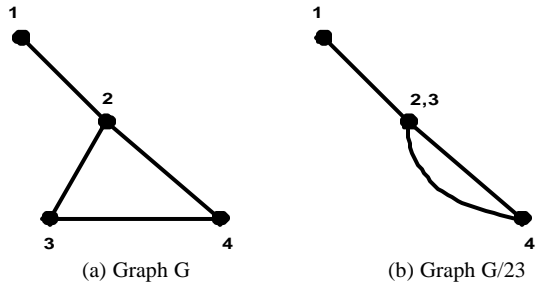

(a) Graph G          (b) Graph G/23

**Fig. 3** - Edge Contraction

## III. RECURSIVE CONTRACTION ALGORITHMS

Given an undirected graph *G(V,E)* we develop an algorithm to scan the set of all cutsets of *G* as $\mathbf{G}=\{\langle S,T \rangle \mid S \subset V, T = V\text{-}S, \langle S \rangle \in C, \langle T \rangle \in C\}$, where *C* is the set of connected graphs.

### A. Simple Partitioning of Vertices

**Algorithm 1.** We divide *V(G)* into two proper disjoint subsets $C_a$ and $C_b$. According to Definition 1, if the subgraphs of *G* that are induced by subsets $C_a$ and $C_b$ are connected subgraphs, the partition $\langle S = C_a, T = C_b \rangle$ is a cutset. In order to find all cutsets, one simple method is to list all possible combinations of vertices that produce two proper disjoint subsets $C_a$ and $C_b$, and then select those that their induced subgraphs $\langle C_a \rangle$ and $\langle C_b \rangle$ are connected [Ahmad 1998; Hongfen 1995]. This simple partitioning algorithm is shown in Fig. 4.

```
Algorithm 1:
a_simple_partitioning_algorithm
{
  input: 1) graph G ;
  output: 1) the list of all cutsets S of G ;
  Initialization: S=∅; partition list F=∅;
  List all possible non-empty partitions of vertices set V(G)
                            in the form of p1 and p2 in F;
  For every partition p1 and p2 in F do
  {
    find induced subgraphs H1=⟨p1⟩ and H2=⟨p2⟩ over G;
    if H1 and H2 both are connected subgraphs
    {
      add partitions p1 and p2 to the cutset list S;
    }
  }
}
```

**Fig. 4** - A simple partitioning Algorithm 1.

**Complexity of Algorithm 1:** The following observation is useful in determining the complexity of Algorithm 1.

**Observation 2.** Dividing *n* objects into two non-empty groups results in $2^{(n-1)} -1$ combinations.

According to Observation 2, one can partition the set of vertices in a graph G with $n$ vertices into $2^{(n-1)}-1$ different combinations. The complexity of checking the connectivity of a graph with $n$ vertices and $m$ edges using depth first search method [Tarjan 1972] is $O(n+m)$. We must check the connectivity of two subgraphs, therefore the complexity of Algorithm 1 is $O((n+m)2^n)$. It grows exponentially with the number of vertices in the graph.

**Lemma 1.** The number of minimal cutsets of a complete graph is $N=2^{(n-1)}-1$.

**Proof:** Each vertex in a complete graph is connected to all other vertices in that graph. Therefore, a partitioning of vertices in a complete graph results in two connected induced subgraphs. According to Observation 2 there are $N=2^{(n-1)}-1$ distinct partitions in $V(G)$. Thus for a complete graph, the number of minimal cutsets is $N=2^{(n-1)}-1$. ÿ

## B. Recursive Contraction Algorithm

**Lemma 2.** A partition $C=\langle S=\{v\}, T = V-\{v\}\rangle$ is a minimal cutset if and only if $v$ is an exterior vertex of $G$.

**Proof:** Let $v$ be an exterior vertex of $G$. According to Definition 2, $H =\langle V-\{v\}\rangle$ is a connected subgraph. The connected subgraphs $F = (\{v\}, Æ)$ and $H$ are two components of $G$ connected to each other via the incident edges in $C$. Therefore, according to Definition 1 the set of edges $C =\langle S=\{v\}, T = V-\{v\}\rangle$ is a minimal cutset for $G$. Conversely, let $C =\langle S=\{v\}, T = V-\{v\}\rangle$ be a minimal cutset of $G$. The subgraph induced by $T_G = V(G)-\{v\}$ must be connected. Therefore, according to Definition 2, $v$ is an exterior vertex of $G$. ÿ

**Observation 3.** A subgraph obtained by deleting two adjacent vertices $u,v\hat{I}\ V$ of a graph $G$ is equivalent to a subgraph obtained by deleting the contracted vertex $g$ from $G/uv$, i.e., $(G/uv)/\{g\}=G/\{u,v\}$. This is a direct result of Definition 4.

**Lemma 3.** For two adjacent vertices $u$ and $v$ in a graph $G$, the partition $\langle S=\{u,v\},T=V-\{u,v\}\rangle$ is a minimal cutset if and only if the contracted vertex $g$ in $G/uv$ is an exterior vertex in $G/uv$.

**Proof:** Let $\langle S=\{u,v\},T =V-\{u,v\}\rangle$ be a minimal cutset for $G$. According to Definition 1, $G/\{u,v\}=\langle V-\{u,v\}\rangle$ must be a connected subgraph. From Observation 3 we know that $(G/uv)/\{g\}=G/\{u,v\}$. This means that the subgraph obtained by deleting the contracted node $g$ is connected, i.e., $g$ is an exterior vertex of $G/uv$. Conversely, let $g$ be an exterior vertex of $G/uv$. This means that $(G/uv)/\{g\}$ is a connected subgraph, i.e., $G/\{u,v\}$ is a connected subgraph. We know that $u$ and $v$ are adjacent vertices, so the subgraph induced by them is connected and contains $u$, $v$, and the edge(s) $uv$. Thus by using Definition 1, we conclude that the partition $\langle S=\{u,v\},T =V(G)-\{u,v\}\rangle$ is a minimal cutset of $G$. ÿ

**Observation 4.** For any subset of vertices $F\tilde{I}\ V$, the induced subgraph $\langle V-F\rangle$ of a given graph $G$ is equivalent to a subgraph obtained by deleting the contracted vertex $g$ from $G/F$, i.e., $(G/F)/\{g\}=\langle V-F\rangle$. This is a direct result of Definition 4.

**Lemma 4.** For any subset of vertices $F\tilde{I}\ V$, the partition $\langle F,V-F\rangle$ is a minimal cutset of graph $G$ if and only if (i) the subgraph induced by $F$ is connected, and (ii) the contracted node $g$ is an exterior vertex of the graph $G/F$.

**Proof:** Let $\langle F,V-F\rangle$ be a minimal cutset for $G$. Then according to Definition 1, the induced subgraphs of $\langle F\rangle$ and $\langle V-F\rangle$ must be connected. From Observation 4, we know that $(G/F)/\{g\} =\langle V-F\rangle$. Thus $(G/F)/\{g\}$ is a connected subgraph. Therefore, according to Definition 2 the contracted vertex $g$ is an exterior vertex in the graph $G/F$. Conversely, let $g$ be an exterior vertex in the graph $G/F$. So the subgraph obtained by deleting the contracted vertex $g$ from $G/F$ (which is equivalent to the induced subgraph $\langle V-F\rangle$) is a connected subgraph. From (i) we know that $\langle F\rangle$ is also a connected subgraph. Therefore, according to Definition 1 we conclude that $\langle F,V-F\rangle$ is a minimal cutset of $G$. ÿ

**Observation 5.** For any subset of vertices $F \subset V$ that has more than one vertex, the subgraph induced by $F$ is connected if and only if there is at least one adjacent vertex for any vertex in $F$. This is a direct result of connectivity and path definitions in graphs.

**Lemma 5.** Finding all connected subgraphs of $G/uv$ that include the contracted vertex $g$ is equivalent to finding all connected subgraphs of $G$ that include $u$ and $v$.

**Proof:** Suppose we find the set of all connected subgraphs $S(g)$ that include the contracted vertex $g$ of the graph $G/uv$, as well as the set of all connected subgraphs $S(u,v)$ that include vertices $u$ and $v$ of the graph $G$. For each connected subgraph $F \in S(g)$ there is a connected subgraph $H \in S(u,v)$ obtained by replacing $g$ with $u$ and $v$ in $F$. Therefore finding $S(g)$ for the graph $G/uv$ is equivalent to finding $S(u,v)$ for $G$. ∎

**Algorithm 2(a).** In this algorithm, we produce a list of all connected subgraphs $S(v)$ in a given graph $G$ that includes vertex $v$. In doing so, we use topological properties of $G$ to reduce the number of partitions whose associated subgraphs must be checked for connectivity. We begin by choosing any vertex of $G$ as *a seed vertex*. We denote all proper subsets of $V$ that include the seed vertex $v$ and induce the connected subgraphs as $S(v)$, i.e., $S(v) = \{ F \subset V \mid v \in F, \langle F \rangle \in C \}$. Suppose the neighborhood of $v$ is $\Gamma(v) = \{u_1, \ldots, u_k\}$. Using Observation 5, we conclude that every proper subset $F \subset V$ in $S(v)$ except $\{v\}$ must include an adjacent vertex $u_i$, i.e.,

$$S(v) = \{v\} \cup S(v,u_1) \cup \ldots \cup S(v,u_k), \quad (1)$$

where $S(v,u_i) = \{F \subset V \mid v \in F, u_i \in \Gamma(v), u_i \in F, F \in C\}$. Now, according to Lemma 5 in order to find $S(v,u_i)$ we must find all connected subgraphs of $G/vu_i$ that include the contracted vertex $g$. But the sets $S(v,u_i)$, $i=1,\ldots,k$ are not disjoint and all partitions will be visited several times. To explain this, suppose that the vertex $z \in \Gamma(v)$ is adjacent to $y \in \Gamma(v)$. Then $\{v,y,z\}$ belongs to both $S(v,z)$ and $S(v,y)$, and all subsets of vertices that include these triple vertices belong to both $S(v,z)$ and $S(v,y)$ as well. To prevent reexamination of recurring subsets, we use the inclusion-exclusion principle and organize vertices according to their position in the BFS tree of the graph. This means that we exclude parts of $S(v,z_i)$ that were previously scanned in $S(v,y)$. To do this we set the following constraint in the scanning algorithm.

**The BFS ordering constraint:** In order to exclude those vertices that have already contributed to the recursion relation (1), we omit all vertices that have a BFS order lower than the BFS order of $v$ in $\Gamma(v)$. Fig. 5 shows the recursive Algorithm 2(a) for finding $S(v)$.

**Theorem 1.** Considering the BFS ordering constraint in Algorithm 2(a), all connected subgraphs of $G$ that include the seed vertex $v$ are scanned once.

**Proof:** Consider a connected subgraph induced by a subset of vertices $F \subset V$, where $v \in F$. There is a unique BFS ordering for all vertices in $F$, denoted as $v$-$u_2$-$\ldots$-$u_k$ where $k<n$. From Observation 5 we know that there is at least one adjacent vertex in $F$ for each vertex in $F$. Algorithm 2(a) scans all adjacent vertices in the vertex list $F$ to construct all possible connected subgraphs in the next recursion level. This ensures that in the recursion level $j$-$1$, a vertex $u_j$ that is adjacent to one of the vertices $v$-$u_2$-$\ldots$-$u_{j-1}$ is added to the list of vertices in $F$, and after $k$ level of recursion the desired connected subgraph is scanned. Since the BFS order of vertices is unique, a connected subgraph cannot be scanned more than once without contradicting the BFS ordering constraint. ∎

**Algorithm 2.** Algorithm 2(a) applies constraint (i) in Lemma 4 and produces a list of all connected subgraphs. To find all cutsets of a graph $G$, we must also apply the second constraint in Lemma 4 and check the connectivity of the induced subgraph $\langle V$-$F \rangle$. The recursive contraction Algorithm 2 depicted in Fig. 6 finds all cutsets in a graph $G$ by generating all possible connected subgraphs $F$ of $G$ and checking that $g$ is an exterior contracted vertex in $G/F$.

```
Algorithm 2(a):
connected_subgraph_recursive_algorithm
{
  inputs: 1) graph G,  2) seed vertex v;
  output: 1) list of all connected subgraphs S(v) of G that include vertex v;
  Initialization: S(v)=Æ, vertex list F={v};
  Find BFS ordering tree of graph G with seed vertex v as its root;
  recursive subroutine: find_all_connected_subgraphs(G, v, F)
  {
    if F is not in S(v)
    {
      add vertex list F to S(v);
    }
    find the neighborhood set ∈(v) for vertex v;
    remove vertices in the vertex list F from ∈(v);
    remove lower BFS ordered vertices of v from ∈(v);
    if ∈(v) is empty
    {
      return;
    }
    else
    {
      recursion loop: for all vertices u of ∈(v) do
      {
        copy the vertex list F to the vertex list H;
        contract the edge (v,u) and find G/uv;
        add vertex u to the  vertex list H;
        set seed vertex v to new contracted vertex g in  G/uv;
        find_all_connected_subgraphs (G/uv, v, H);
      }
    }
  }
}
```

**Fig. 5** - Algorithm 2(a) for listing all connected subgraphs that include the seed vertex $v$

**Theorem 2.** Considering the BFS ordering constraint in Algorithm 2, all cutsets of a given graph $G$ are scanned once.

**Proof:** Each cutset divides a given graph $G$ into two unique connected subgraphs, one of which includes the seed vertex $v$. In order to specify a cutset, it is sufficient to find its corresponding connected subgraph that includes the seed vertex $v$. We know from Theorem 1 that Algorithm 2(a) scans all connected subgraphs once, which include the seed vertex $v$. Since Algorithm 2 uses the same approach and examines all possible connected subgraphs, we conclude that this algorithm scans all cutsets of the graph $G$ once. ÿ

**Complexity of Algorithm 2:** Algorithm 2 discards disconnected subgraphs that cannot generate any cutset. In order to determine the complexity of this algorithm, we enumerate the number of disconnected subgraphs of $G$ and subtract it from the number of all subgraphs of $G$ that was computed in Lemma 2. We use Lemma 6 to do this.

**Lemma 6.** The number of disconnected subgraphs of a graph $G$ that include the seed vertex $v$ and exclude all neighbors of $v$ is $2^{c-1}$, where $c = n-|A(v)|-1$, $n$ is the total number of vertices in $G$ and $A(v)$ is the number of adjacent vertices of $v$.

**Proof:** Let $(S,T)$ be a partition of the vertex set $V$. Now, $v \hat{I} S$ and all vertices $u_i$ in $A(v)$ belong to $T$. Then $c = n-A(v)-1$ remaining vertices can be divided between $S$ and $T$ in an arbitrary manner. Since the number of non-empty subsets of the set of remaining vertices is $2^{c-1}$, then the number of disconnected subgraphs induced by $S$ is also $2^{c-1}$. ÿ

```
Algorithm 2:
all_cutset_recursive_algorithm
{
  inputs: 1) graph G,  2) seed vertex v;
  output: 1) list of all cutsets S(v) of G that include vertex v;
  Initialization: S(v)=Æ; vertex list F={v};
  Find BFS ordering tree of the graph G with seed vertex v as its root;
  recursive subroutine: find_all_cutsets(G,v,F)
  {
    if vertex v is an exterior vertex of the graph G
    {
      if F is not in S(v)
      {
        add the vertex list F to S(v);
      }
    }
    find the neighborhood set Γ(v) for vertex v;
    remove vertices in the vertex list F from Γ(v);
    remove lower BFS ordered vertices of v from Γ(v);
    if Γ(v) is empty
    {
      return;
    }
    else
    {
      recursion loop: for all vertices u of Γ(v) do
      {
        copy the vertex list F to the vertex list H;
        contract the edge (v,u) and find G/uv;
        add vertex u to the vertex list H;
        set the seed vertex v to the new contracted vertex g of G/uv;
        find_all_cutsets(G/uv,v,H);
      }
    }
  }
}
```

**Fig. 6** - The recursive contraction Algorithm 2 for listing all cutsets of a given graph *G*.

**Algorithm 2(b).** Lemma 6 computes the number of disconnected subgraphs that have not been scanned in each recursion of Algorithm 2. In order to find the total number of disconnected subgraphs, we accumulate this value in every recursion. Algorithm 2(b) shown in Fig. 7 gives the number of disconnected subgraphs.

The number of disconnected subgraphs depends on the structure of the BFS tree of the seed vertex. Lemma 6 indicated that the seed vertex with the lowest number of neighbors results in the lowest complexity of Algorithm 2.

## C. Enhanced Recursive Contraction Algorithm

Algorithm 2 is adopted form of Tsukiyama's algorithm and its complexity is linear per cutset. In this section, we develop Algorithm 3, the complexity of which is also linear per cutset, but is less than that of Algorithm 2. To do this, we first examine the notion of 1-point extension of a cutset introduced by Tsukiyama *et al*. [1980].

```
Algorithm 2(b):
 Enumerating_disconnected_subgraph_algorithm
{
 inputs: 1) graph G,  2) seed vertex v;
 output: 1) disn number of all disconnected Subgraphs of G that include the seed vertex v ;
 Initialization: disn=0, vertex list F={v};
  n= |V(G)|;
 Find BFS ordering tree of the graph G with seed vertex v as its root;
 recursive subroutine: enumerate_partial_disconnected_subgraphs(G,v,F)
 {
   find the neighborhood set G(v) for vertex v;
   remove vertices in vertex list F from G(v) ;
   remove lower BFS ordered vertices of v from G(v) ;
   if G(v) is empty
   {
     return;
   }
   else
   {
    c= n-|G (v)| -1;
    disn=disn+2^{c-1};
   recursion loop: for all vertices u of G(v) do
   {
    copy the vertex list F to vertex list H;
    contract the edge (v,u) and find G/uv;
    add vertex u to vertex list H;
    set seed vertex v to the new contracted vertex g of G/uv;
    enumerate_partial_disconnected_subgraphs(G/uv,v,H);
   }
  }
 }
}
```

**Fig. 7**- Algorithm 2(b) for enumerating disconnected subgraphs that include the seed vertex *v*

**Definition 5. The 1-point extension:** Let $X \subseteq Y \subseteq V$. If both $\langle X,X'\rangle$ and $\langle Y,Y'\rangle$ are cutsets, then there exist a $v \in Y-X$ so that $\langle X+v, X'-v\rangle$ is a cutset [Whited *et al.* 1990]. The set $X+v$ is called a 1-point extension of $X$.

According to the above definition if we find a method to determine the extension nodes of a given cutset $\langle X,X'\rangle$, then we can construct other cutsets such as $\langle Y,Y'\rangle$ by transferring one of these extension nodes from the sub-partition $X$ to the sub-partition $X'$. Algorithm 2 contains this method. In this algorithm we begin from any seed vertex and find all 1-point extension by examining all adjacent vertices. The crucial step in Algorithm 2 is to determine the set $W(X)$ comprised of all vertices $v$ for which $X+v$ is a 1-point extension of $X$, i.e., $W(X) = \{ v \notin X' / v \in G(X), v \in K(X')\}$, where $K(X')$ is the set of all exterior vertices of $\langle X'\rangle$. It is important to note that determining $K(X')$, and thus $W(X)$, for each cutset $\langle X,X'\rangle$ is the most computation-intensive part of Algorithm 2. The main idea in our enhanced algorithm is to develop an efficient way to determine whether or not a given vertex is in $K(X')$. Prior to this, we define the following terms.

**Definition 6. A cluster and the pivot:** According to Definition 2, when we delete an interior vertex from a graph $G$, more than one distinct component remains in the resulting graphs. Each one of such components is called a cluster, and the deleted interior vertex is called the pivot vertex of these clusters.

**Definition 7. The absorbable cluster:** A cluster in which all its vertices have a BFS order greater than the BFS order of its pivot (the BFS order of the last vertex contracted in the pivot) is an absorbable cluster. Otherwise it is an inabsorbable cluster (Fig. 8).

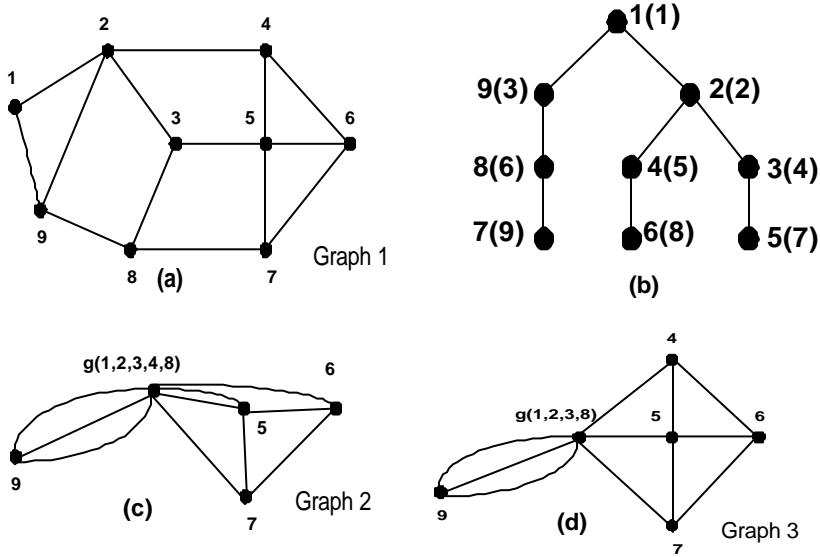

**Fig. 8** - (a) A sample graph, (b) The breadth first search ordering *m* of node *n* is shown in the form of *n(m)*, (c) Graph 2 is generated by contracting nodes 1,2,3,4,8. The BFS order of the contracted node g(1,2,3,4,8) is 6 (BFS order of the last node, i.e., node 8). The cluster <9> is inabsorbable but the cluster <5,6,7> is absorbable, (d) Graph 3 is generated by contracting nodes 1,2,3,8. The BFS order of the contracted node is 6. Both clusters <9> and <4,5,6,7> are inabsorbable, because the BFS order of node 4 is 5, which is lower than the BFS order of the pivot.

**Theorem 3.** A 1-point extension for a pivot vertex with more than one inabsorbable cluster does not exist.

**Proof:** Contraction of vertices in Algorithm 2 eventually results in complete absorption of all absorbable clusters in a contraction node. When only one cluster remains, the contracted vertex is transformed into an exterior vertex of the graph. But, if there is more than one inabsorbable cluster, contraction of vertices in Algorithm 2 does not transform the contracted vertex into an exterior vertex of the graph. This is due to the fact that there is at least one vertex in each inabsorbable cluster whose BFS order is lower than the BFS order of the last vertex contracted in the pivot vertex. According to the BFS ordering constraint, this vertex cannot be selected for contraction, resulting in at least two distinct clusters in the graph. Therefore, the pivot vertex always remains an interior vertex of the graph, and consequently it cannot belong to *W(X)*. ÿ

**Observation 9.** A pivot vertex with one inabsorbable cluster is transformed into an exterior vertex when the vertex with the highest BFS order among all vertices of absorbable clusters is contracted in the pivot vertex.

**Algorithm 3:** As shown in Fig. 9, we use Theorem 3 to modify Algorithm 2 to prevent consecutive contraction when the contracted vertex is a pivot vertex with more than one inabsorbable clusters. When the contracted vertex is a pivot vertex with one inabsorbable cluster, we use Observation 9; and this algorithm goes through a cut-through transient phase, until the contracted vertex transforms into an exterior vertex. At this time, the algorithm returns to a normal recursive form. During the cut-through transient phase, the contraction process is carried out as before, but the recursion, which is the most time consuming part of the algorithm and checks the exteriority of the contracted node, is not performed. This is due to the fact that according to Observation 9, the contraction node is an interior node. The end of this phase is detected by the contraction of the highest BFS order vertex among all vertices of the absorbable clusters.

```
Algorithm 3:
Enhanced_all_cutset_recursive_algorithm
{
 inputs: 1) graph G,  2)seed vertex v;
 output: 1) list of all cutsets S(v) of G that include the vertex v;
 Initialization: S(v)=Æ; vertex list F={v}; dummy_flag=FALSE;
 Find BFS ordering tree of the graph G with seed vertex v as its root;
 recursive subroutine: find_all_cut-sets(G,v,F)
 {
  if vertex v is an exterior vertex of the graph G
  {
   if F is not in S(v)
   {
    add the vertex list F to S(v);
   }
  }
   else
   {
    find the set of clusters clst for the pivot vertex v;
    if there is more than one inabsorbable cluster in clst
    {
     return;
    }
    else
    {
     if there is one inabsorbable cluster in clst
     {
      highest_BFS_order_vertex = the highest BFS order vertex among
                     vertices of absorbable clusters;
      dummy_flag = TRUE;
     return;
     }
    }
   }
  find the neighborhood set Γ(v) for the vertex v;
  remove vertices in the vertex list F from Γ(v);
  remove the lower BFS ordered vertices of v from Γ(v);
  if Γ(v) is empty
  {
   return;
  }
   else
   {
    recursion loop: for all vertices u of Γ(v) do
    {
     copy the vertex list F to the vertex list H;
     contract the edge (v,u) and find G/uv;
     add the vertex u to the vertex list H;
     set the seed vertex v to new contracted vertex g of G/uv;
     if (dummy_flag is FALSE )
     {
      find_all_cutsets(G/uv,v,H);
     }
      else
      {
       if (v = the highest_BFS_order_vertex) dummy_flag = FALSE;
      }
    }
   }
  }
 }
}
```

**Fig. 9** - The enhanced recursive contraction Algorithm 3 for scanning all minimal cutsets of a given graph *G*

**Complexity of Algorithm 3:** Algorithm 3, in addition to discarding disconnected subgraphs, also discards a class of connected subgraphs that according to Theorem 3 and Observation 9 cannot specify a cutest of the graph. In order to determine the complexity of Algorithm 3, we enumerate the number of connected subgraphs of $G$ that are ignored in this algorithm. We use lemma 7 and 8 to do this.

**Lemma 7.** Suppose that during the execution of Algorithm 3, a pivot with more than one inabsorbable cluster is obtained. Assume that there are $j < n$ vertices previously contracted in the pivot and the highest BFS order among them is $m \geq j$. The number of ignored connected subgraphs for this pivot is in the order of $O(2^{n-m})$.

**Proof:** We have a connected subgraph $H$ with $n-j+1$ vertices, among them $m-j$ vertices have a BFS order lower than $m$ and only $n-m$ vertices contract in a pivot (the remaining node) due to BFS ordering constraint. We ignore all connected subgraphs of $H$ that are generated by absorbing these $n-m$ vertices. From the complexity of Algorithm 2 we observe that the number of connected subgraphs for a given connected graph is $k*2^{n-1}$, where $n$ is number of vertices of the graph and $k < 1$ is a constant that depends on the selection of seed vertex and the topology of the graph. Therefore, we conclude that the number of discarded connected subgraphs in Algorithm 3 is in the order of $O(2^{n-m})$. ÿ

**Lemma 8.** Suppose that during the execution of Algorithm 3, a pivot with exactly one inabsorbable cluster is obtained. Assume that there are $j < n$ vertices previously contracted in the pivot and the highest BFS order among them is $m \geq j$. Also assume that the highest BFS order among the vertices of the absorbable clusters is $f \leq n$. The number of ignored connected subgraphs for this pivot is in the order of $O(2^{f-m})$.

**Proof:** We have a connected subgraph $H$ with $n-j+1$ vertices, among them $m-j$ vertices have a BFS order lower than $m$ and cannot be absorbed in a pivot according to the BFS ordering constraint. Also, $n-f+1$ vertices have a BFS order greater than or equal to $f$ that can generate the 1-point extensions. Therefore, only $f-m$ vertices can be absorbed in the pivot during the cut-through transient phase. We discard all connected subgraphs of $H$ that are obtained from the absorption of these vertices. From the complexity of Algorithm 2 we know that the number of connected subgraphs obtained from a combination of $n$ vertices is in the order of $O(2^{n-1})$. Therefore, we conclude that the number of connected subgraphs that are ignored in Algorithm 3 when we encounter a pivot with exactly one inabsorbable cluster is in the order of $O(2^{f-m})$. ÿ

The complexity of Algorithm 3 is lower than that of Algorithm 2. This is due to the fact that the search space in Algorithm 3 does not include several groups of connected subgraphs that cannot generate any new cutsets. The number of such groups increases exponentially with the size of the graph.

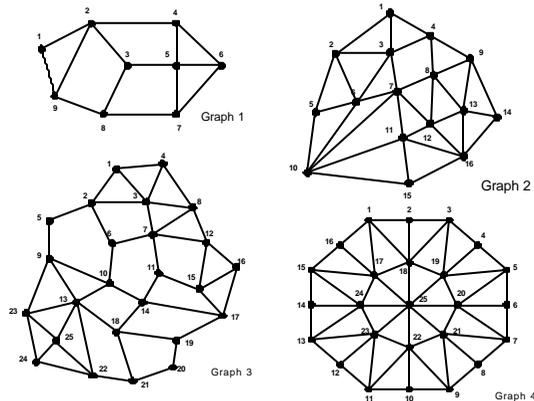

**Fig. 10** - Some sample graphs for performance comparison of the algorithms

## IV. PERFORMANCE COMPARISON

In this Section, we present the results of implementing the simple Algorithm 1, the recursive contraction Algorithm 2, and its enhanced version (Algorithm 3). In order to compare the complexities of these algorithms we test them for a number of sample graphs as depicted in Fig. 10. The number of cutsets for Graphs 1,2,3 and 4 are 66, 2232, 17518 and 28448 respectively.

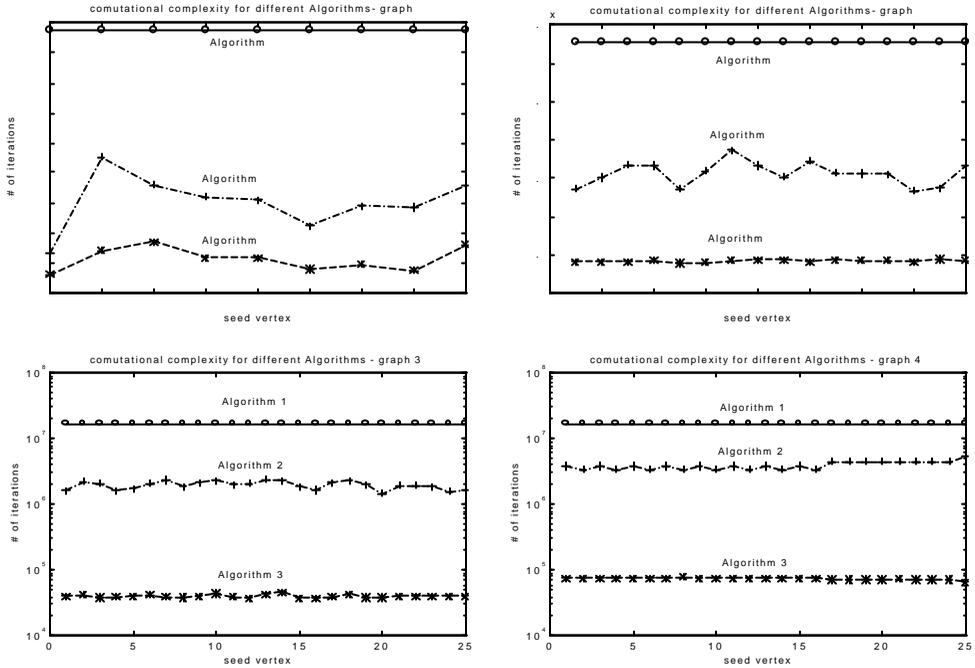

**Fig. 11**- Comparison of complexities of Algorithms 1,2 and 3 for sample graphs in Fig. 10

All three algorithms give us the correct values of cutsets irrespective of the seed vertex. We calculate the number of iterations in for different seed vertices. The complexities of Algorithms 1, 2 and 3 are compared in Fig. 11. It is evident that the complexity of Algorithm 3 is substantially lower compared to Algorithms 1 and 2 for all seed vertices. The complexity of Algorithms 2 and 3 for larger graphs is reduced several orders of magnitude compared to Algorithm 1. In Algorithm 2, the minimum number of iterations is obtained by choosing the lowest degree node as seed vertex. This is not the case for Algorithm 3, for which the variance of the number of iterations as a function of the seed vertex in Algorithm 3 is low compared to Algorithm 2.

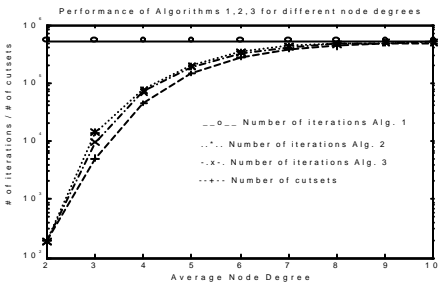

**Fig. 12** - The number of iterations and cutests vs. average node degree for a network of size 20.

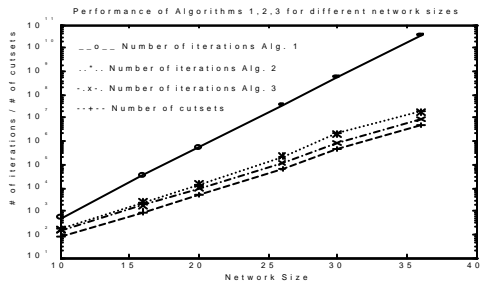

**Fig. 13 -** The number of iterations and cutsets vs. average network size for a node degree of size 3.

To show that the number of iterations in Algorithms 2 and 3 is proportionate to the number of cutsets in a graph, we generate random graphs with different sizes and average node degrees, and apply Algorithms 2 and 3 to them. Fig. 12 shows the results for a graph with 20 nodes and different average node degrees from 2 to 10. In Fig. 13 we show the results for random graphs with an average node degree equal to 3 and different sizes from 10 to 36 nodes. For each sample point in Figs. 12 and 13, we apply Algorithms 2 and 3 ten times for different random monolithic graphs and calculate the mean values for the iteration numbers as well as the number of cutsets. In applying Algorithm 2, we begin with finding the least degree node of the graph and use this node as the seed vertex. It is evident from Figs. 12 and 13 that for Algorithms 2 and 3, the number of iterations follows the number of cutsets. We conclude that the complexities of both Algorithms 2 and 3 are proportionate to the number of cutsets irrespective of the size of the graph but the complexity of Algorithm 3 is less than the complexity of Algorithm 2. By increasing the node degree the difference in the performances of Algorithms 2 and 3 becomes smaller and for larger node degrees both asymptotically tend towards the performance of Algorithm 1.

Figs. 14 and 15 show the ratio of mean number of iterations to the mean number of cutsets, which we call the efficiency factor. As this value gets closer to 1, the efficiency of the algorithm is improved. Figs. 14 and 15 show that Algorithm 3 is more efficient for all node degree values and network sizes. As shown in Fig. 15, the efficiency factor of Algorithm 3 is independent of graph size, and is about 1.8 when the average node degree is 3, but the efficiency factor of Algorithm 2 decreases for larger network sizes. The difference in the performances of Algorithms 2 and 3 is more significant for large bounded degree mesh networks, which represents the majority of telecommunications networks.

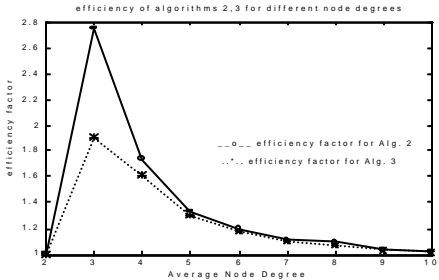
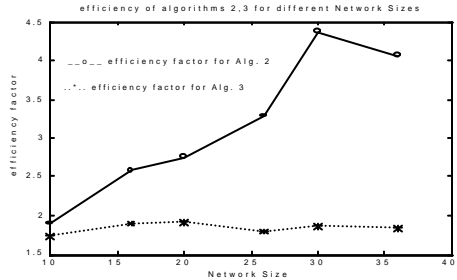

**Fig. 14** – The efficiency factor for the network of Fig. 12 vs. the average node degree

**Fig. 15** – The efficiency factor for the network of Fig. 13 vs. network size.

## V. CONCLUSIONS

In this paper we presented a pseudo-polynomial recursive algorithm for scanning all minimal cutsets of a given graph. This algorithm is based on consecutive contraction of edges in the graph and uses the concept of interior vertices for checking the connectivity of the induced subgraphs that are generated from partitioning of the graph. We used the BFS ordering of vertices to prevent revisiting a possible state more than once. We also used the concept of an absorbable cluster for an interior vertex (pivot) in a graph to substantially reduce the complexity of the algorithm.

We applied our proposed algorithms to different graphs and showed that scanning all cutsets of these graphs is done in linear time per cutset. Since the number of minimal cutsets increases exponentially with network size (Fig. 13), we conclude that the complexity of Algorithm 3 is pseudo-polynomial. We also applied Algorithm 3 to monolithic randomly generated graphs with different sizes and average node degrees and showed that the efficiency factor improves for graphs with large average node degrees, but does not change considerably with network size when the average node degree is kept constant.

## REFERENCES

AHMAD, S.H. 1988. Simple enumeration of minimal cutsets of acyclic directed graph. *IEEE Trans. On Reliability*, *37 (5),* 484-487.

AHUJA, R.R., MAGNATI, T.L., AND ORLIN, J.B. 1993. *Network Flows- Theory, Algorithms, and Applications,* Prentice-Hall International.

BALL, M.O. (ed.) 1995. *Network Models (chapter on network reliability).* Elsevier, Amsterdam.

COLBOURN, C.J. 1987. *The Combinatorics of Network Reliability. vol. 4 of The International Series of Monographs on Computer Science*, Oxford University Press.

FARD, N.S., AND LEE, T.H. 1999. Cutset Enumeration of network systems with link and node failure. *Reliability Engineering and System Safety, 65,* 141-146.

GOLDBERG, A.V., AND TARJAN, R.E. 1988. A new approach to the maximum flow problem. *Journal of the ACM, 35,* 921-940.

HONGFEN, Z. 1995. A Simple enumeration of all minimal cutsets of an all-terminal graph. *The Journal of China Universities of Posts and Telecommunications, 2(2).*

KARGER, D.R. 2000. Minimum cuts in near-linear time. *Journal of the ACM*, *47,* 46-76.

PROVAN, J.S., AND BALL, M.O. 1983. The complexity of counting cuts and of computing the probability that a graph is connected. *SIAM Journal of Computing*, *12,* 777-788.

PROVAN, J.S., AND SHIER, D.R. 1996. A paradigm for listing (s,t)-cuts in graphs. *Algorithmica, 15,* 351-372.

RAI, S. 1982. A cutset approach to reliability evaluation in communication networks. *IEEE Trans. on Reliability, 31,* 428-431.

SHARAFAT, A.R., AND MA'ROUZI, O.R. 2002. The most congested cut-set: deriving a lower bound for the chromatic number in the RWA problem. Submitted for publication to the *IEEE Communications Letters*.

SHIER, D.R. 1988. Algebraic Aspects of Computing Network Reliability. In *Applications of Discrete Mathematics,* R.D. Ringeisen and F.S. Roberts (eds.), SIAM, Philadelphia, pp. 135-147.

TARJAN, R.E. 1972. Depth first search and linear graph algorithms. *SIAM Journal of Computing*, *1,* 146-160.

TSUKIYAMA, S., SHIRAKAWA, I., OZAKI, H., AND ARIYOSHI, H. 1980. An algorithm to enumerate all cutsets of a graph in linear time per cutset *Journal of the ACM*, *27,* 619-632.

WHITED, D.E., SHIER, D.R., AND JARVIS, J.P. 1990. Reliability computations for planar networks. *OSRA Journal of Computing*, *2(1):* 46-60.